\documentclass[11pt]{article}
\usepackage{graphicx}
\usepackage{color}
\usepackage{indentfirst}
\usepackage{amsmath,amssymb,amsfonts,amsthm,bm}
\usepackage{cases}
\usepackage{cite}
\usepackage{mathrsfs}

\newtheorem{theorem}{Theorem}[section]
\newtheorem{definition}{Definition}[section]

\newtheorem{corollary}{Corollary}[section]
\newtheorem{example}{Example}[section]
\newtheorem{remark}{Remark}[section]
\input{tcilatex}

\begin{document}

\title{Best proximity point results for almost contraction and application to nonlinear differential equation}
\author{Azhar Hussain \thanks{Department of Mathematics, University of Sargodha, Sargodha-40100, Pakistan. Email: hafiziqbal30@yahoo.com,}
Mujahid Abbas \thanks{Department of Mathematics, Government College University, Lahore 54000, Pakistan. Email: abbas.mujahid@gmail.com},
Muhammad Adeel\thanks{Department of Mathematics, University of Sargodha, Sargodha-40100, Pakistan. Email: adeel.uosmaths@gmail.com,} Tanzeela Kanwal\thanks{Department of Mathematics, University of Sargodha, Sargodha-40100, Pakistan. Email: tanzeelakanwal16@gmail.com}.}
\date{}
\maketitle

\begin{center}
\noindent\textbf{\ Abstract}
\end{center}
Brinde [Approximating fixed points of weak
contractions using the Picard itration, Nonlinear Anal. Forum {\bf 9}
(2004), 43-53] introduced almost contraction mappings and proved Banach
contraction principle for such mappings. The aim of this paper is to
introduce the notion of multivalued almost $\Theta$-contraction mappings and
present some best proximity point results for this new class of mappings. As
applications, best proximity point and fixed point results for weak single
valued $\Theta$-contraction mappings are obtained. An example is presented
to support the results presented herein. An application to a nonlinear
differential equation is also provided.
 \newline



\section{Introduction and preliminaries}

The following concept was introduced by Berinde as `weak contraction' in
\cite{BR}. But in \cite{BR1}, Berinde renamed `weak contraction' as `almost
contraction' which is appropriate.

\begin{definition} Let $(X,d)$ be a metric space. A mapping $%
F:X\rightarrow X$ is called almost contraction or $(\delta ,L)$-contraction
if there exist a constant $\delta \in (0,1)$ and some $L\geq 0$ such that
for any $x,y\in X,$ we have
\begin{equation}
d(Fx,Fy)\leq \delta .d(x,y)+Ld(y,Fx),
\end{equation}
\end{definition}
Von Neumann \cite{von} considered fixed points of multivalued
mappings in the study of game theory. Indeed, the fixed point results for
multivalued mappings play a significant role in study of control theory and
in solving many problems of economics and game theory.

Nadler \cite{nadler} used the concept of the Hausdorff metric to obtain
fixed points of multivalued contraction mappings and obtained the Banach
contraction principle as a special case.\newline
Here, we recall that a Hausdorff metric $H$ induced by a metric $d$ on a set
$X$ is given by%
\begin{equation}
H(A,B)=\max \{\sup_{x\in A}d(x,B),\sup_{y\in B}d(y,A)\}
\end{equation}%
for every $A,B\in CB(X)$, where $CB(X)$ is the collection of the closed and
bounded subsets of $X.$

M. Berinde and V. Berinde \cite{Berinde2} introduced the notion of
multivalued almost contraction as follows:

Let $(X,d)$ be a metric space. A mapping $F:X\rightarrow CB(X)$ is called
almost contraction if there exist two constants $\delta \in (0,1)$ and $%
L\geq 0$ such that for any $x,y\in X,$ we have%
\begin{equation}
H(Fx,Fy)\leq \delta d(x,y)+LD(y,Fx).
\end{equation}

Berinde \cite{Berinde2} proved Nadler's fixed point theorem in (\cite{nadler}%
):

\begin{theorem} Let $(X,d)$ be a complete metric space and $%
F:X\rightarrow CB(X)$ a almost contraction. Then $F$ has a fixed point..
\end{theorem}
Jleli {\it et al.} \cite{jleli} defined $\Theta $-contraction mapping as
follows:

A mapping $F:X\rightarrow X$ is called $\Theta$-contraction if for any $%
x,y\in X$
\begin{equation}
\Theta (d(Fx,Fy))\leq \lbrack \Theta (d(x,y))]^{k}  \label{E1.1}
\end{equation}%
where, $k\in (0,1)$ and $\Theta :(0,\infty )\rightarrow (1,\infty )$ is a
mapping which satisfies the following conditions

\begin{description}
\item[($\Theta_{1}$)] $\Theta$ is nondecreasing;

\item[($\Theta _{2}$)] for each sequence $\{\alpha _{n}\}\subseteq {\Bbb R}%
^{+}$, $\lim\limits_{n\rightarrow \infty }\Theta (\alpha _{n})=1$ if and
only if $\lim\limits_{n\rightarrow \infty }(\alpha _{n})=0$;

\item[($\Theta _{3}$)] there exist $0<k<1$ and $l\in (0,\infty )$ such that $%
\lim\limits_{\alpha \rightarrow 0^{+}}\frac{\Theta (\alpha )-1}{\alpha ^{k}}%
=l$;
\end{description}

Denote%
\begin{equation}
\Omega =\{\Theta :(0,\infty )\rightarrow (1,\infty ):\Theta ~~ satisfies~~
\Theta _{1}-\Theta _{3}\}.
\end{equation}

\begin{theorem} (\cite{jleli}) Let $(X,d)$ be a complete metric
space and $F:X\rightarrow X$ a $\Theta $-contraction, then $F$ has a unique
fixed point.
\end{theorem}
Hancer {\it et al.} \cite{Hancer} introduced the notion of multi-valued $%
\Theta$-contraction mapping as follows:

Let $(X,d)$ be a complete metric space and $T:X\rightarrow CB(X)$ a
multivalued mapping. Suppose that there exists $\Theta \in \Omega$ and $%
0< k <1$ such that
\begin{equation}
\Theta (H(Tx,Ty))\leq \lbrack \Theta (d(x,y))]^{k}
\end{equation}%
for any $x,y$ $\in $ $X$ provided that $H(Tx,Ty)>0$, where $CB(X)$ is a
collection of all nonempty closed and bounded subsets of $X$.

\begin{theorem} \label{multi} Let $(X,d)$ be a complete metric
space and $F:X\rightarrow K(X)$ a multi-valued $\Theta$-contraction, then $F$
has a fixed point.
\end{theorem}
Let $A$ and $B$ be two nonempty subsets of a metric space $(X,d)$ and $%
F:A\rightarrow CB(B)$. A point $x^{\ast }\in A$ is called a best proximity
point of $F$ if
\begin{eqnarray*}
D(x^{\ast }, Fx^{\ast }) &=&{\it \inf }\{d(x^{\ast },y):y\in Fx^{\ast
}\}=dist(A,B),
\end{eqnarray*}%
where
\[
dist(A,B) =\inf \{d(a,b):a\in A,b\in B\}.
\]
If $A\cap B\neq \phi ,$ then $x^{\ast }$ is a fixed point of $F.$ If $A\cap
B=\phi $, then $D(x,Fx)>0$ for all $x\in A$ and $F$ has no fixed point.

Consider the following optimization problem:%
\begin{equation}
{\it \min }\{D(x,Fx):x\in A\}.
\end{equation}%
It is then important to study necessary conditions so that the above
minimization problem has at least one solution.

Since
\begin{equation}
d(A,B)\leq D(x,Fx)
\end{equation}%
for all $x\in A$. Hence the optimal solution to the problem
\begin{equation}
{\it \min }\{D(x,Fx):x\in A\}
\end{equation}%
for which the value $d(A,B)$ is attained is indeed a best proximity point of
multivalued mapping $F.$

For more results in this direction, we refer to \cite{Shahzad, Gabeleh1,
Gabeleh, Amini, Basha, Veeramani, Choudhury, Dimri, Eldred, Eldred1, azhar, LAH,
Plebaniak, Plebaniak1} and references mentioned therein.

Let $A$ and $B$ two nonempty subsets of $X.$ Denote
\begin{eqnarray*}
A_{0} &=&\{a\in A:d(a,b)=d(A,B)~~ \text{for some} ~~b\in B\} \\
B_{0} &=&\{b\in B:d(a,b)=d(A,B)~~\text{for some} ~~a\in A\}.
\end{eqnarray*}
\begin{definition} \cite{KSankar} Let $(X,d)$ be a metric space
and $A_{0}\neq \phi $, we say that the pair $(A,B)$ has the $P$-property if%
\begin{equation}
\left.%
\begin{array}{c}
d(x_{1},y_{1})=d(A,B) \\
d(x_{2},y_{2})=d(A,B)%
\end{array}
\right\} ~~ implies ~ that ~~d(x_{1},x_{2})=d(y_{1},y_{2}),
\end{equation}
where $x_{1},x_{2}\in ~A$ and $y_{1},y_{2}\in ~B$$.$
\end{definition}
\begin{definition}
\cite{zhng} Let $(X,d)$ be a metric space and $%
A_{0}\neq \phi $, we say that the pair $(A,B)$ has the weak $P$-property if%
\begin{equation}
\left.%
\begin{array}{c}
d(x_{1},y_{1})=d(A,B) \\
d(x_{2},y_{2})=d(A,B)%
\end{array}%
\right\}~~ implies~ that ~~d(x_{1},x_{2})\leq d(y_{1},y_{2}),
\end{equation}%
where $x_{1},x_{2}\in ~A$ and $y_{1},y_{2}\in ~B.$
\end{definition}
\begin{definition}
\cite{BSch} Let $(X,d)$ be a metric space, $%
A,B $ two subsets of $X$ and $\alpha :A\times A\rightarrow \lbrack 0,\infty
) $. A mapping $F:A\rightarrow 2^{B}\backslash \{\phi \}$ is called $\alpha$%
-proximal admissible if%
\begin{equation}
\left.%
\begin{array}{c}
\alpha (x_{1},x_{2})\geq 1, \\
d(u_{1},y_{1})=d(A,B), \\
d(u_{2},y_{2})=d(A,B)%
\end{array}%
\right\} ~~implies~ that ~~\alpha (u_{1},u_{2})\geq 1
\end{equation}%
for all $x_{1},x_{2},u_{1},u_{2}\in A,~y_{1}\in Fx_{1}$ and $y_{2}\in Fx_{2}$.
\end{definition}
\begin{definition}\cite{BSch}
Let $F:X\rightarrow CB(Y)$ be a
multi-valued mapping, where $(X,d_{1})$, $(Y,d_{2})$ are two metric spaces.
A mapping $F$ is said to be continuous at $x\in X$ if $H(Fx,Fx_{n})%
\rightarrow 0$ whenever $d_{1}(x,x_{n})\rightarrow 0$ as $n\rightarrow
\infty .$
\end{definition}
The aim of this paper is to obtain some best proximity point results for multivalued almost
$\Theta$-contraction mappings. We also present some best proximity point and fixed point results for single valued mappings. Moreover, an example to prove the validity and application to nonlinear differential equation for the usability of our results is presented. Our results extend, unify and generalize the
comparable results in the literature.
\section{Best proximity points of multivalued almost $\Theta $%
-contraction}

We begin with the following definition:

\begin{definition}\label{def1}
Let $A,B$ be two nonempty
subsets of a metric space $(X,d)$ and $\alpha :A\times A\rightarrow \lbrack
0,\infty )$. Let $\Theta \in \Omega $ be a continuous function. A
multivalued mapping $F:A\rightarrow 2^{B}\backslash \{\phi \}$ is called
almost $\Theta $-contraction if for any $x,y\in A,$ we have%
\begin{equation}
\alpha (x,y)\Theta \lbrack H(Fx,Fy)]\leq \lbrack \Theta (d(x,y)+\lambda
D(y,Fx)))]^{k}  \label{Eq100}
\end{equation}%
where $k\in (0,1)$ and $\lambda \geq 0.$
\end{definition}
\begin{theorem} \label{Th1} Let $(X,d)$ be a complete metric
space and $A$, $B$ nonempty closed subsets of $X$ such that $A_{0}\neq \phi $. Suppose that $F:A\rightarrow K(B)$ is a continuous mapping such that
\begin{description}
\item[(i)] $Fx\subseteq B_{0}$ for each $x\in A_{0}$ and $(A,B)$ satisfies
the weak $P$-property;

\item[(ii)] $F$ is $\alpha$-proximal admissible mapping;

\item[(iii)] there exists $x_{0},x_{1}\in A_{0}$ and $y_{0}\in
Fx_{0}\subseteq B_{0}$ such that \newline
$d(x_{1},y_{0})=d(A,B)$ and $\alpha (x_{0},x_{1})\geq 1$;

\item[(iv)] $F$ is multivalued almost $\Theta $-contraction.
\end{description}
Then $F$ has a best proximity point in $A$.
\begin{proof} Let $x_{0},x_{1}$ be two given points in $A_{0}$
and $y_{0}\in Fx_{0}\subseteq B_{0}$ such that $d(x_{1},y_{0})=d(A,B)$ and $%
\alpha (x_{0},x_{1})\geq 1.$ If $y_{0}\in Fx_{1}$, then $d(A,B)\leq
D(x_{1},Fx_{1})\leq d(x_{1},y_{0})=d(A,B)$ implies that $%
D(x_{1},Fx_{1})=d(A,B)$ and $x_{1}$ is a best proximity point of $F$. If $%
y_{0}\notin Fx_{1}$ then,
\[
0<D(y_{0},Fx_{1})\leq H(Fx_{0},Fx_{1}).
\]%
Since $F(x_{1})\in K(B),$ we can choose $y_{1}\in Fx_{1}$ such that
\[
1<\Theta \lbrack d(y_{0},y_{1})]\leq \Theta \lbrack H(Fx_{0},Fx_{1})].
\]%
As $F$ is multivalued almost $\Theta $-contraction mapping, we have
\begin{eqnarray}
1 &<&\Theta \lbrack d(y_{0},y_{1})]\leq \alpha (x_{0},x_{1})\Theta \lbrack
H(Fx_{0},Fx_{1})]  \nonumber  \label{Eq101} \\
&\leq &[\Theta (d(x_{0},x_{1})+\lambda D(x_{1},Fx_{0}))]^{k}  \nonumber \\
&=&[\Theta (d(x_{0},x_{1}))]^{k}.
\end{eqnarray}%
Since $y_{1}\in Fx_{1}\subseteq B_{0}$, there exists $x_{2}\in A_{0}$ such
that $d(x_{2},y_{1})=d(A,B)$ and $\alpha (x_{1},x_{2})\geq 1$. By weak $P$%
-property of the pair $(A,B)$ we obtain that $d(x_{2},x_{1})\leq
d(y_{0},y_{1})$. If $y_{1}\in Fx_{2},$ then $x_{2}$ is a best proximity
point of $F$. If $y_{1}\notin Fx_{2}$, then
\[
D(y_{1},Fx_{2})\leq H(Fx_{1},Fx_{2}).
\]%
We now choose $y_{2}\in Fx_{2}$ such that
\begin{eqnarray}
1<\Theta \lbrack d(y_{1},y_{2})] &\leq &\Theta \lbrack H(Fx_{1},Fx_{2})]
\nonumber  \label{Eq102} \\
&\leq &\alpha (x_{1},x_{2})\Theta \lbrack H(Fx_{1},Fx_{2})]  \nonumber \\
&\leq &[\Theta (d(x_{1},x_{2})+\lambda D(x_{2},Fx_{1}))]^{k}  \nonumber \\
&=&[\Theta (d(x_{1},x_{2}))]^{k}.
\end{eqnarray}%
Continuing this process, we can obtain two sequences $\{x_{n}\}$ and $%
\{y_{n}\}$ in $A_{0}\subseteq A$ and $B_{0}\subseteq B,$ respectively such
that $y_{n}\in Fx_{n}$ and it satisfies
\[
d(x_{n+1},y_{n})=d(A,B)~~\text{with}~~\alpha (x_{n},x_{n+1})\geq 1
\]%
where $n=0,1,2,....$ Also,
\begin{eqnarray}
1<\Theta \lbrack d(y_{n},y_{n+1})] &\leq &\Theta \lbrack H(Fx_{n},Fx_{n+1})]
\nonumber  \label{Eq103} \\
&\leq &\alpha (x_{n},x_{n+1})\Theta \lbrack H(Fx_{n},Fx_{n+1})]  \nonumber \\
&\leq &[\Theta (d(x_{n},x_{n+1})+\lambda D(x_{n},Fx_{n+1}))]^{k}  \nonumber
\\
&=&[\Theta (d(x_{n},x_{n+1}))]^{k}.
\end{eqnarray}%
implies that
\begin{equation}  \label{Eq1005}
1<\Theta \lbrack d(y_{n},y_{n+1})]\leq (\Theta (d(x_{n},x_{n+1})))^{k}.
\end{equation}
Since
\begin{equation}
d(x_{n+1},y_{n})=d(A,B)  \label{Eq104}
\end{equation}%
and
\begin{equation}
d(x_{n},y_{n-1})=d(A,B)
\end{equation}%
for all $n\geq 1$, it follows by the weak $P$-property of the pair $(A,B)$
that
\begin{equation}  \label{Eq105}
d(x_{n},x_{n+1})\leq d(y_{n-1},y_{n})
\end{equation}
for all $n\in {\Bbb N}$. \ Now by repeated application of (\ref{Eq1005}), (%
\ref{Eq105}) and the monotone property of $\Theta $, we have
\begin{eqnarray}
1<\Theta \lbrack d(x_{n},x_{n+1})] &\leq &\Theta (d(y_{n-1},y_{n}))\leq
\Theta (H(Fx_{n-1},Fx_{n}))  \nonumber  \label{Eq35} \\
&\leq &\alpha (x_{n-1},x_{n})\Theta (H(Fx_{n-1},Fx_{n}))  \nonumber \\
&\leq &[\Theta (d(x_{n-1},x_{n})+\lambda D(x_{n},Fx_{n-1}))]^{k}  \nonumber
\\
&=&(\Theta (d(x_{n-1},x_{n})))^{k}\leq (\Theta (d(y_{n-2},y_{n-1})))^{k}
\nonumber \\
&\leq &(\Theta (H(Fx_{n-2},Fx_{n-1})))^{k}  \nonumber \\
&\leq &(\alpha (x_{n-2},x_{n-1})\Theta (H(Fx_{n-2},Fx_{n-1}))^{k}  \nonumber
\\
&\leq &[\Theta (d(x_{n-2},x_{n-1})+\lambda D(x_{n-1},Fx_{n-2}))]^{k^{2}}
\nonumber \\
&=&(\Theta (d(x_{n-2},x_{n-1})))^{k^{2}}  \nonumber \\
&&.  \nonumber \\
&&.  \nonumber \\
&&.  \nonumber \\
&\leq &(\Theta (d(x_{0},x_{1})))^{k^{n}}.
\end{eqnarray}%
for all $n\in {\Bbb N}\cup \{0\}$. This shows that $\lim\limits_{n%
\rightarrow \infty }\Theta (d(x_{n},x_{n+1}))=1$ and hence $%
\lim_{n\rightarrow \infty }d(x_{n},x_{n+1})=0.$ From ($\Theta _{3}$), there
exist $0<r<1$ and $0<l\leq \infty $ such that
\begin{equation}
\lim_{n\rightarrow \infty }\frac{\Theta (d(x_{n},x_{n+1}))-1}{%
[d(x_{n},x_{n+1})]^{r}}=l.  \label{ch3eq3}
\end{equation}%
Assume that $l<\infty $ and $\beta =l/2$. From the definition of the limit
there exists $n_{0}\in {\Bbb N}$ such that
\[
\left\vert \frac{\Theta (d(x_{n},x_{n+1}))-1}{[d(x_{n},x_{n+1})]^{r}}%
-l\right\vert \leq B,\text{~for~all~}n\geq n_{0}
\]%
which implies that
\[
\frac{\Theta (d(x_{n},x_{n+1})-1}{[d(x_{n},x_{n+1})]^{r}}\geq l-\beta =\beta
\text{ for all }n\geq n_{0}.
\]%
Hence
\[
n[d(x_{n},x_{n+1})]^{r}\leq n\alpha \lbrack \Theta (d(x_{n},x_{n+1})-1]\text{
for all }n\geq n_{0},
\]%
where $\alpha =1/\beta $. \ Assume that $l=\infty $. \ Let $\beta >0$ be a
given real number. From the definition of the limit there exists $n_{0}\in
{\Bbb N}$ such that
\[
\frac{\Theta (d(x_{n},x_{n+1})-1}{[d(x_{n},x_{n+1})]^{r}}\geq \beta \text{
for~all~~}n\geq n_{0}
\]
implies that
\[
n[d(x_{n},x_{n+1})]^{r}\leq n\alpha \lbrack \Theta (d(x_{n},x_{n+1})-1]\text{
for~all~}~n\geq n_{0},
\]%
where $\alpha =1/\beta $. Hence, in all cases there exist $\alpha >0$ and $%
n_{0}\in {\Bbb N}$ such that
\[
n[d(x_{n},x_{n+1})]^{r}\leq n\alpha \lbrack \Theta (d(x_{n},x_{n+1})-1]\text{%
~for~all~~}n\geq n_{0}.
\]%
From (\ref{Eq35}), we have
\[
n[d(x_{n},x_{n+1})]^{r}\leq n\alpha \lbrack \Theta (d(x_{n},x_{n+1})-1]\text{
}~\text{for~all~~}n\geq n_{0}.
\]%
On taking the limit as $n\rightarrow \infty $ on both sides of the above
inequality, we have
\begin{equation}
\lim_{n\rightarrow \infty }n[d(x_{n},x_{n+1})]^{r}=0.  \label{ch3eq4}
\end{equation}%
It follows from (\ref{ch3eq4}) that there exists $n_{1}\in {\Bbb N}$ such
that
\[
n[d(x_{n},x_{n+1})]^{r}\leq 1\text{ for~all~~}n>n_{1}.
\]%
This implies that
\[
d(x_{n},x_{n+1})\leq \frac{1}{n^{1/r}}~\text{\ for~all~}~n>n_{1}.
\]%
Now, for $m>n>n_{1}$, we have
\[
d(x_{n},x_{m})\leq \sum_{i=n}^{m-1}d(x_{i},x_{i+1})\leq \sum_{i=n}^{m-1}%
\frac{1}{i^{\frac{1}{r}}}.
\]%
Since $0<r<1,$ $\sum_{i=n}^{\infty }\frac{1}{i^{\frac{1}{r}}}$ converges.
Therefore $d(x_{n},x_{m})\rightarrow 0$ as $m,n\rightarrow \infty $.\newline
This shows that $\{x_{n}\}$ and $\{y_{n}\}$ are Cauchy sequences in $A$ and $%
B,$ respectively. Next, we assume that there exists elements $u,v\in A$ such
that
\[
x_{n}\rightarrow u~~\text{and~}~y_{n}\rightarrow v~~as~~n\rightarrow \infty
.
\]%
Taking limit as $n\rightarrow \infty $ in (\ref{Eq104}), we obtain that
\begin{equation}  \label{Eq108}
d(u,v)=d(A,B).
\end{equation}
Now , we claim that $v\in Tu$. Since $y_{n}\in Fx_{n}$, we have
\[
D(y_{n},Fu)\leq H(Fx_{n},Fu).
\]%
Taking limit as $n\rightarrow \infty $ on both sides above sides of above
inequality, we have\newline
\[
D(v,Fu)=\lim_{n\rightarrow \infty }D(y_{n},Fu)\leq \lim_{n\rightarrow \infty
}H(Fx_{n},Fu)=0.
\]%
As $Fu\in K(B)$, $D(v,Fu)=0$ implies that $v\in Fu$. \ By (\ref{Eq108}), we
have
\[
D(u,Fu)\leq d(u,v)=dist(A,B)\leq D(u,Fu),
\]%
which implies that $D(u,Fu)=dist(A,B)$ and hence $u$ is a best proximity
point of $F$ in $A$.
\end{proof}
\end{theorem}
\begin{remark} In the next theorem, we replace the continuity
assumption on $F$ with the following condition:

If $\{x_{n}\}$ is a sequence in $A$ such that $\alpha (x_{n},x_{n+1})\geq 1$
for all $n$ and $x_{n}\rightarrow x\in A$ as $n\rightarrow \infty $, then
there exists a subsequence $\{x_{n(k)}\}$ of $\{x_{n}\}$ such that $\alpha
(x_{n(k)},x)\geq 1$ for all $k$. If the above condition is satisfied then we
say that the set $A$ satisfies $\alpha -$ subsequential property.
\end{remark}
\begin{theorem} \label{Th2} Let $(X,d)$ be a complete metric
space and $(A,B)$ a pair of nonempty closed subsets of $X$ such that $%
A_{0}\neq \phi $. \ Let $F:A\rightarrow K(B)$ be a multivalued mapping such
that conditions (i)-(iv) of Theorem \ref{Th1} are satisfied. Then $F$ has
a best proximity point in $A$ provided that $A$ satisfies $\alpha -$
subsequential property.
\begin{proof} Following Arguments similar to those in the proof
of Theorem \ref{Th1}, we obtain two sequences $\{x_{n}\}$ and $\{y_{n}\}$ in
$A$ and $B,$ respectively such that
\begin{equation}
\alpha (x_{n},x_{n+1})\geq 1,
\end{equation}%
\begin{equation}
x_{n}\rightarrow u\in A~~\text{and~}~y_{n}\rightarrow v\in B~~\text{as}%
~~n\rightarrow \infty
\end{equation}%
and
\begin{equation}
d(u,v)=dist(A,B)
\end{equation}%
By given assumption, there exists a subsequence $\{x_{n(k)}\}$ of $\{x_{n}\}$
such that $\alpha (x_{n(k)},u)\geq 1$ for all $k$. \ Since $y_{n(k)}\in
Fx_{n(k)}$ for all $k\geq 1$, applying condition (iv) of Theorem \ref{Th1}%
, we obtain that
\begin{eqnarray}
1<\Theta (D(y_{n(k)},Fu)) &\leq &\Theta (H(Fx_{n(k)},Fu))  \nonumber
\label{Eq109} \\
&\leq &\alpha (x_{n(k)},u)\Theta (H(Fx_{n(k)},Fu))  \nonumber \\
&=&(\Theta (d(x_{n(k)},u)))^{k}.
\end{eqnarray}%
On taking limit as $k\rightarrow \infty $ in (\ref{Eq109}) and using the
continuity of $\Theta $, we have $\Theta (D(v,Fu)=1$. Therefore, by ($\Theta
_{2}$) we obtain that $D(v,Fu)=0$. As shown in the proof of Theorem \ref{Th1}%
, we have $D(v,Fu)=dist(A,B)$ and hence $u$ is a best proximity point of $F$
in $A$.
\end{proof}
\end{theorem}
\begin{remark} \label{R6} To obtain the uniqueness of the best
proximity point of multivalued almost $\Theta $-contraction mappings, we
propose the following ${\cal H}$ condition:

${\cal H}:$ for any best proximity points $x_{1},x_{2}$ of mapping $F,$ we
have
\[
\alpha (x_{1},x_{2})\geq 1.
\]
\end{remark}
\begin{theorem} \label{Th3} Let $A$ and $B$ be two nonempty
closed subsets of a complete metric space $(X,d)$ such that $A_{0}\neq \phi $
and $F:A\rightarrow K(B)$ continuous multivalued mapping satisfying the
conditions of Theorem \ref{Th1} (respectively in Theorem \ref{Th2}). Then the
mapping $F$ has a unique best proximity point provided that it satisfies the
condition ${\cal H}$.

\begin{proof} Let $x_{1},x_{2}$ be two best proximity points of $%
F $ such that $x_{1}\neq x_{2}$, then by the given hypothesis ${\cal H}$ we
have $\alpha (x_{1},x_{2})\geq 1$ and $%
D(x_{1},Fx_{1})=D(A,B)=D(x_{2},Fx_{2}) $. Since $Fx_{1}$ and $Fx_{2}$ are
compact sets, there exists elements $y_{0}\in Fx_{1}$ and $y_{1}\in Fx_{2}$
such that
\[
d(x_{1},y_{0})=dist(A,B)\text{ and }d(x_{2},y_{1})=dist(A,B).
\]%
As $F$ satisfies the weak $P$-property, we have
\[
d(x_{1},x_{2})\leq d(y_{0},y_{1}),
\]%
Since $F$ is multivalued almost $\Theta $-contraction mapping, we obtain
that
\begin{eqnarray*}
\Theta (d(x_{1},x_{2}))\leq \Theta (d(y_{0},y_{1})) &\leq &\Theta
(H(Fx_{1},Fx_{2})) \\
&\leq &\alpha (x_{1},x_{2})\Theta (H(Fx_{1},Fx_{2})) \\
&\leq &[\Theta (d(x_{1},x_{2})+\lambda D(x_{2},Fx_{1}))]^{k} \\
&=&(\Theta (d(x_{1},x_{2})))^{k} \\
&<&\Theta (d(x_{1},x_{2})),
\end{eqnarray*}%
a contradiction. Hence $d(x_{1},x_{2})=0$, and $x_{1}=x_{2}.$
\end{proof}
\end{theorem}
If the pair $(A,B)$ satisfies the weak $P$-property, then it satisfies the $%
P $-property, we have the following corollaries:

\begin{corollary} \label{c1} Let $(X,d)$ be a complete metric
space and $(A,B)$ a pair of nonempty closed subsets of $X$ such that $%
A_{0}\neq \phi $. Suppose that a continuous mapping $F:A\rightarrow K(B)$
satisfies the following properties:

\begin{description}
\item[(i)] $Fx\subseteq B_{0}$ for each $x\in A_{0}$ and $(A,B)$ satisfies
the $P$-property;

\item[(ii)] $F$ is multivalued $\alpha$-proximal admissible mapping;

\item[(iii)] there exists $x_{0},x_{1}\in A_{0}$ and $y_{0}\in
Fx_{0}\subseteq B_{0}$ such that $d(x_{1},y_{0})=d(A,B)$ and $\alpha
(x_{0},x_{1})\geq 1$;

\item[(iv)] $F$ is multivalued almost $\Theta $-contraction.
\end{description}

Then $F$ has a best proximity point in $A$.
\end{corollary}
\begin{corollary}
 \label{c2} Let $(X,d)$ be a complete metric
space and $(A,B)$ a pair of nonempty closed subsets of $X$ such that $%
A_{0}\neq \phi $. let $F:A\rightarrow K(B)$ be a multi-valued mapping such
that conditions (i)-(iv) of Corollary \ref{c1} are satisfied. Then $F$ has
a best proximity point in $A$ provided that $A$ has $\alpha -$ subsequential
property.
\end{corollary}
Now we give an example to support Theorem \ref{Th1}.

\begin{example} \label{exam1} \ Let $X={\Bbb R}^{2}$ be a usual
metric space. Let
\begin{equation}
A=\{(-2,2),(2,2),(0,4)\}
\end{equation}
and
\begin{equation}
B=\{(-8,\gamma ):\gamma \in \lbrack -8,0]\}\cup \{(8,\gamma ):\gamma \in
\lbrack -8,0]\}\cup \{(\beta ,-8):\beta \in (-8,8)\}.
\end{equation}
Then $d(A,B)=8$, $\ A_{0}=\{(-2,2),(2,2)\}$ and $B_{0}=\{(-8,0),(8,0)\}$.

Define the mapping $F:A\rightarrow K(B)$ by
\[
Fx=\left\{%
\begin{array}{cc}
\{(-8,0)\} & if ~~x=(-2,2) \\
\{(8,0)\} & if ~~x=(2,2) \\
\{(\beta ,-8):\beta \in (-8,8)\} & if~~x=(0,4). \\
&
\end{array}%
\right.
\]
and $\alpha :A\times A\rightarrow \lbrack 0,\infty )$ by
\begin{equation}
\alpha ((x,y),(u,v))=\frac{11}{10}.
\end{equation}%
Clearly, $F(A_{0})\subseteq B_{0}.$ For $(-2,2),(2,2)\in A$ and $%
(-8,0),(8,0)\in B,$ we have%
\[
\left\{%
\begin{array}{c}
d((-2,2),(-8,0))=d(A,B)=8, \\
d((2,2),(8,0))=d(A,B)=8.%
\end{array}%
\right.
\]
Note that
\begin{equation}
d((-2,2),(2,2))<d((-8,0),(8,0)).
\end{equation}%
that is, the pair $(A,B)$ has weak $P$-property. Also, $F$ is $\alpha $%
-proximal admissible mapping. Now we show that $F$ is multivalued almost $%
\Theta $-contraction where $\Theta :(0,\infty )\rightarrow (1,\infty )$ is
given by $\Theta (t)=5^{t}$.

Note that
\begin{equation}
\alpha ((-2,2),(2,2))\Theta \lbrack H(F(-2,2),F(2,2))]=\frac{11}{10}(5^{16})
\label{ExE1}
\end{equation}%
and
\begin{equation}
\lbrack \Theta (d((-2,2),(2,2))+\lambda D((2,2),(-8,0)))]^{k}=\left( \
5^{28}\right) ^{k}.  \label{ExE2}
\end{equation}%
If we take $k\in (\frac{717}{1250},1)$ and $\lambda =2$ in $(\ref{ExE2})$,
we have
\begin{equation}
\frac{11}{10}(5^{16})<\left( \ 5^{28}\right) ^{k}.
\end{equation}%
Similarly,%
\begin{equation}
\alpha (x,y)\Theta \lbrack H(Fx,Fy)]\leq \lbrack \Theta (d(x,y)+\lambda
D(y,Fx)))]^{k}
\end{equation}
holds for the remaining pairs. Hence all the conditions of Theorem \ref{Th1}
are satisfied. Moreover, $(-2,2)$, $(2,2)$ are best proximity points of $F$
in $A$.
\end{example}
\begin{remark} Note that mapping $F$ in the above example does
not hold for the case of Nadler \cite{nadler} as well as for Hancer {\it et al.}
\cite{Hancer}. For if, take $x=(-2,2),y=(2,2)\in A$, we have
\[
\Theta (H(Tx,Ty))=5^{16}>5^{4}>(5^{4})^{k}=[\Theta (d(x,y))]^{k}
\]%
for $k\in (0,1)$. Also
\[
H(Tx,Ty)=16>4=d(x,y)>\alpha d(x,y)
\]%
for $\alpha \in (0,1)$.
\end{remark}
\begin{remark}
 In above example \ref{exam1}, the pair $(A,B) $
does not satisfy the $P$-property and hence the Corollary \ref{c1} is not
applicable in this case.
\end{remark}

\section{Application to single valued mappings}

In this section, we obtain some best proximity point results for
singlevalued mappings as applications of our results obtain in section 2.

\begin{definition} \cite{jleli} \ Let $(X,d)$ be a metric space
and $A,B$ two subsets of $X$, a nonself mapping $F:A\rightarrow B$ is called
$\alpha $-proximal admissible if
\begin{equation}
\begin{cases} \alpha(x_{1},x_{2})\geq 1, \\
d(u_{1},Fx_{1})=d(A,B),~~~\text{implies}~~~\alpha(u_{1},u_{2})\geq 1 \\
d(u_{2},Fx_{2})=d(A,B).\end{cases}
\end{equation}%
for all $x_{1},x_{2},u_{1},u_{2}\in A$ where $\alpha :A\times A\rightarrow
\lbrack 0,\infty )$.
\end{definition}
\begin{definition}
 Let $\alpha :A\times A\rightarrow \lbrack
0,\infty )$ and $\Theta :(0,\infty )\rightarrow (1,\infty )$ a nondecreasing
and continuous function. A mapping $F:A\rightarrow B$ is called almost $%
\Theta $-contraction if for any $x,y\in A,$ we have
\begin{equation}
\alpha (x,y)\Theta (d(Fx,Fy))\leq \lbrack \Theta (d(x,y)+\lambda
d(y,Fx)))]^{k}  \label{Eq1001}
\end{equation}%
where, $k\in (0,1)$ and $\lambda \geq 0.$
\end{definition}
\begin{theorem} \label{L1.1} Let $(X,d)$ be a complete metric
space and $(A,B)$ a pair of nonempty closed subsets of $X$ such that $A_{0}$
is nonempty. If $F:A\rightarrow B$ is a continuous mapping such that

\begin{description}
\item[(i)] $F(A_{0}) \subseteq B_{0}$ and $(A, B)$ satisfies the weak $P$%
-property;

\item[(ii)] $F$ is $\alpha$-proximal admissible mapping;

\item[(iii)] there exists $x_{0}, x_{1}\in A_{0}$ such that $d(x_{1},
Fx_{0})=d(A, B)$ and $\alpha(x_{0}, x_{1}) \geq 1$;

\item[(iv)] $F$ is almost $\Theta$-contraction.
\end{description}

Then $F$ has a best proximity point in $A$.

\begin{proof} As for every $x\in X,\{x\}$ is compact in $X$.
Define a multivalued mapping $T:A\rightarrow K(B)$ by $Tx=\{Fx\}$ for $x\in
A $. The continuity of $F$ implies that $T$ is continuous. Now $%
F(A_{0})\subseteq B_{0}$ implies that $Tx=\{Fx\}\subseteq B_{0}$ for each $%
x\in A_{0}.$ If $x_{1},x_{2},v_{1},v_{2}\in A,$ $y_{1}\in Tx_{1}=\{Fx_{1}\}$
and $y_{2}\in Tx_{2}=\{Fx_{2}\}$ are such that\newline
\begin{equation}
\alpha (x_{1},x_{2})\geq 1,~~d(v_{1},y_{1})=dist(A,B)~~\text{and}%
~~d(v_{2},y_{2})=dist(A,B).
\end{equation}%
That is,
\begin{equation}
\alpha (x_{1},x_{2})\geq
1,~~d(v_{1},Fx_{1})=dist(A,B)~~and~~d(v_{2},Fx_{2})=dist(A,B).
\end{equation}%
Then we have $\alpha (v_{1},v_{2})\geq 1$ as $F$ is $\alpha $-proximal
admissible mapping. Hence $T$ is $\alpha $-proximal admissible mapping.

Suppose there exist $x_{0},x_{1}\in A_{0}$ such that $%
d(x_{1},Fx_{0})=dist(A,B)$ and $\alpha (x_{0},x_{1})\geq 1$. Let $y_{0}\in
Tx_{0}=\{Fx_{0}\}\subseteq B_{0}$. Then $d(x_{1},Fx_{0})=dist(A,B)$ gives
that $d(x_{1},y_{0})=dist(A,B)$. By condition (iii), there exist $%
x_{0},x_{1}\in A_{0}$ and $y_{0}\in Tx_{0}\subseteq B_{0}$ such that $%
d(x_{1},y_{0})=dist(A,B)$ and $\alpha (x_{0},x_{1})\geq 1.$

Since $F$ is almost $\Theta $-contraction, we have
\begin{equation}
\alpha (x,y)\Theta (H(Tx,Ty))=\alpha (x,y)\Theta \lbrack d(Fx,Fy)]\leq
\lbrack \Theta (d(x,y)+\lambda D(y,Fx)))]^{k},
\end{equation}%
for any $x,y\in A$ which implies that $T$ is multivalued almost $\Theta $%
-contraction. Thus, all the conditions of Theorem \ref{Th1} are satisfied
and hence $T$ has a best proximity point $x^{\ast }$ in $A$. Thus we have $%
D(x^{\ast },Tx^{\ast })=dist(A,B)$ and hence $d(x^{\ast },Fx^{\ast
})=dist(A,B)$, that is $x^{\ast }$ is a best proximity point of $F$ in $A$.
\end{proof}
\end{theorem}
\begin{theorem} \label{L1.2} Let $(X,d)$ be a complete metric
space and $(A,B)$ a pair of nonempty closed subsets of $X$ such that $A_{0}$
is nonempty. If $F:A\rightarrow B$ is a mapping such that conditions
(i)-(iv) of Theorem $\ref{L1.1}$ are satisfied. Then $F$ has a best
proximity point in $A$ provided that $A$ satisfies $\alpha$-subsequential
property.

\begin{proof} Let $T:A\rightarrow K(B)$ be as given in proof of
Theorem \ref{L1.1}. Following arguments similar to those in the proof of
Theorem \ref{L1.1}, we obtain that

\begin{description}
\item[(i)] $Tx \subseteq B_{0}$ for each $x_{0} \in A_0$;

\item[(ii)] $T$ is multi-valued $\alpha$-proximal admissible mapping;

\item[(iii)] there exists $x_{0}, x_{1}\in A_{0}$ and $y_0 \in Tx_0
\subseteq B_0$ such that $d(x_1, y_0)=d(A, B)$ and $\alpha(x_{0}, x_{1})
\geq 1$;

\item[(iv)] $T$ is multivalued almost $\Theta $-contraction.
\end{description}

Thus, all the conditions of Theorem \ref{Th2} are satisfied and hence $T$
has a best proximity point $x^{\ast }$ in $A,$ that is,
\[
D(x^{\ast },Tx^{\ast })=dist(A,B)
\]%
Consequently, $d(x^{\ast },Fx^{\ast })=dist(A,B)$, and $x^{\ast }$ is a best
proximity point of $F$ in $A$.
\end{proof}
\end{theorem}
\begin{corollary} \label{c3} Let $(X,d)$ be a complete metric
space and $(A,B)$ a pair of nonempty closed subsets of $X$ such that $A_{0}$
is nonempty. If $F:A\rightarrow B$ is a continuous mapping such that

\begin{description}
\item[(i)] $F(A_{0}) \subseteq B_{0}$ and $(A, B)$ satisfies the $P$%
-property;

\item[(ii)] $F$ is $\alpha$-proximal admissible mapping;

\item[(iii)] there exists $x_{0}, x_{1}\in A_{0}$ such that $d(x_{1},
Tx_{0})=d(A, B)$ and $\alpha(x_{0}, x_{1}) \geq 1$;

\item[(iv)] $F$ is almost $\Theta$-contraction.
\end{description}

Then $F$ has a best proximity point in $A$.

\begin{proof} Replace the condition of weak P-property with
P-property in Theorem \ref{L1.1}.
\end{proof}
\end{corollary}
\begin{corollary} \label{c4} Let $(X,d)$ be a complete metric
space and $(A,B)$ be a pair of nonempty closed subsets of $X$ such that $%
A_{0}$ is nonempty. If $F:A\rightarrow B$ is a mapping such that conditions
(i)-(iv) of Corollary $\ref{c3}$ are satisfied. Then $F$ has a best
proximity point in $A$ provided that $A$ satisfies $\alpha -$ subsequential
property.
\begin{proof} Replace the condition of weak P-property with
P-property in Theorem \ref{L1.2}.
\end{proof}
\end{corollary}
\section{Fixed point results for single and multi-valued mappings}
In this section, fixed points of singlevalued and multivalued almost $\Theta
$-contraction mappings are obtained.

Taking $A=B=X$ in Theorem \ref{Th1} (Theorem \ref{Th2}), we obtain
corresponding fixed point results for almost $\Theta $-contraction mappings.

\begin{theorem} \label{Th5.1} Let $(X,d)$ be a complete metric
space. If $F:X\rightarrow K(X)$ is a continuous mapping satisfying

\begin{description}
\item[(i)] $F$ is $\alpha$ admissible mapping;

\item[(ii)] there exists $x_{0}\in X$ such that $\alpha(x_{0}, Fx_{0})\geq 1$%
;

\item[(iii)] $F$ is multivalued almost $\Theta $-contraction.
\end{description}

Then $F$ has a fixed point in $X$.
\end{theorem}
\begin{theorem} \label{Th5.2} Let $(X,d)$ be a complete metric
space. Let $F:X\rightarrow K(X)$ be a multi-valued mapping such that
conditions (i)-(iii) of Theorem \ref{Th5.1} are satisfied. Then $F$ has a
fixed point in $X$ provided that $X$ satisfies $\alpha -$ subsequential
property.
\end{theorem}
Taking $A=B=X$ in Theorem \ref{L1.1} (in Theorem \ref{L1.2}), we obtain the
corresponding fixed point results of almost $\Theta $-contraction mappings.

\begin{theorem}\label{Th5.3} Let $(X,d)$ be a complete metric
space. Let $F:X\rightarrow X$ be a continuous mapping satisfying

\begin{description}
\item[(i)] $F$ is $\alpha$ admissible mapping;

\item[(ii)] there exists $x_{0}\in X$ such that $\alpha(x_{0}, Fx_{0})\geq 1$%
;

\item[(iii)] $F$ is almost $\Theta$-contraction.
\end{description}

Then $F$ has a fixed point in $X$.
\end{theorem}
\begin{theorem} \label{Th5.4} Let $(X,d)$ be a complete metric
space. Let $F:X\rightarrow X$ be a multi-valued mapping such that conditions
(i)-(iii) of Theorem \ref{Th5.1} are satisfied. Then $F$ has a fixed point
in $X$ provided that $X$ has a $\alpha -$subsequential property.
\end{theorem}
\begin{remark} In Theorem 4.1 (respectively in 4.3)

(i). If we take $\alpha(x, y)=1$, we obtain the main results of Durmaz \cite%
{durmaz} and Altun \cite{altun}.

(ii). Taking $\lambda=0$ and $\alpha(x, y)=1$, we obtain the main result of
Hancer {\it et al.} \cite{Hancer} and Jelli \cite{jleli} respectively.

(iii). Taking $\alpha(x, y)=1$ and $\Theta(t)=e^{t}$, we obtain the main
result of Berinde \cite{Berinde2} and \cite{BR}.

(iv). Taking $\alpha(x, y)=1$, $\lambda=0$ and $\Theta(t)=e^{t}$, we obtain
the main result of Nadler \cite{nadler} and Banach \cite{banach}.
\end{remark}
\section{Application to Nonlinear Differential Equations}

Let $C([0,1])$ be the set of all continuous functions defined on $[0,1]$ and
$d:C([0,1])\times C([0,1])\rightarrow {\Bbb R}$ be the metric defined by
\begin{equation}
d(x,y)=||x-y||_{\infty }=\max_{t\in \lbrack 0,1]}|x(t)-y(t)|.  \label{1}
\end{equation}%
It is known that $(C([0,1]),d)$ is a complete metric space.

Let us consider the two-point boundary value problem of the second-order
differential equation:
\begin{equation}
\left.
\begin{array}{cc}
-\frac{d^{2}x}{dt^{2}}=f(t,x(t)) & t\in \lbrack 0,1]; \\
x(0)=x(1)=0 &
\end{array}%
\right\}  \label{2}
\end{equation}%
where $f:[0,1]\times {\Bbb R}\rightarrow {\Bbb R}$ is a continuous mapping.

The Green function associated with (\ref{2}) is defined by
\begin{equation}
G(t,s)=\left\{
\begin{array}{cc}
t(1-s) & if~~0\leq t\leq s\leq 1, \\
s(1-t) & if~~0\leq s\leq t\leq 1.%
\end{array}%
\right.
\end{equation}%
Let $\phi :{\Bbb R}\times {\Bbb R}\rightarrow {\Bbb R}$ be a given function.

Assume that the following conditions hold:

\begin{enumerate}
\item[(i)] $|f(t,a)-f(t,b)|\leq \max\limits_{a, b\in {\Bbb R}}|a-b|$ for all
$t\in \lbrack 0,1]$ and $a, b\in {\Bbb R}$ with $\phi (a,b)\geq 0$;

\item[(ii)] {\it there exists }$x_{0}\in C[0,1]$ such that $\phi
(x_{0}(t),Fx_{0}(t))\geq 0${\it \ for all }$t\in \lbrack 0,1]$ where $F:C[0,1]\rightarrow C[0,1]$;

\item[(iii)] {\it for each }$t\in \lbrack 0,1]${\it \ and }$x,y\in C[0,1]$%
{\it , }$\phi (x(t),y(t))\geq 0${\it \ implies }$\phi (Fx(t), Fy(t))\geq 0$%
{\it ; }

\item[(iv)] {\it for each }$t\in \lbrack 0,1]${\it , if }$\{x_{n}\}${\it \
is a sequence in }$C[0,1]${\it \ such that }$x_{n}\rightarrow x${\it \ in }$%
C[0,1]${\it \ and }$\phi (x_{n}(t),x_{n+1}(t))\geq 0${\it \ for all }$n\in N$%
{\it , then }$\phi (x_{n}(t),x(t))\geq 0${\it \ for all }$n\in N${\it . }
\end{enumerate}

We now prove the existence of a solution of the second order differential
equation (\ref{2}).

\begin{theorem} Under the assumptions (i)-(iv), (\ref{2}) has a
solution in $C^2([0, 1])$.
\begin{proof} It is well known that $x\in C^{2}([0,1])$ is a
solution of (\ref{2}) is equivalent to $x\in C([0,1])$ is a solution of the
integral equation
\begin{equation}
x(t)=\int_{0}^{1}G(t,s)f(s,x(s))ds,~t\in \lbrack 0,1].  \label{3}
\end{equation}%
Let $F:C[0,1]\rightarrow C[0,1]$ be a mapping defined by
\begin{equation}
Fx(t)=\int_{0}^{1}G(t,s)f(s,x(s))ds.  \label{4}
\end{equation}%
Suppose that $x,y\in C([0,1])$ such that $\phi (x(t),y(t))\geq 0$ for all $%
t\in \lbrack 0,1]$. By applying (i), we obtain that
\begin{eqnarray*}
|Fu(x)-Fv(x)| &=&\int_{0}^{1}G(t,s)f(s,x(s))ds-\int_{0}^{1}G(t,s)f(s,y(s))ds
\\
&\leq &\int_{0}^{1}G(t,s)[f(s,x(s))-f(s,y(s))]ds \\
&\leq &\int_{0}^{1}G(t,s)|f(s,x(s))-f(s,y(s))|ds \\
&\leq &\int_{0}^{1}G(t,s)\cdot (\max |x(s)-y(s)|)ds \\
&\leq &||x-y||_{\infty }\cdot \sup_{t\in \lbrack 0,1]}\Bigg(%
\int_{0}^{1}G(t,s)ds\Bigg).
\end{eqnarray*}%
Since $\int\limits_{0}^{1}G(t,s)ds=-(t^{2}/2)+(t/2)$, for all $t\in \lbrack
0,1]$, we have
\[
\sup_{t\in \lbrack 0,1]}\Bigg(\int_{0}^{1}G(t,s)ds\Bigg)=\frac{1}{8}.
\]%
It follows that
\begin{equation}
||Fx-Fy||_{\infty }\leq \frac{1}{8}(||x-y||_{\infty }  \label{5}
\end{equation}%
Taking exponential on the both sides, we have
\begin{eqnarray}
e^{||Fx-Fy||_{\infty }} &\leq &e^{\frac{1}{8}(||x-y||_{\infty })}  \nonumber
\label{6} \\
&=&[e^{(||x-y||_{\infty })}]^{\frac{1}{8}},
\end{eqnarray}%
for all $x,y\in C[0,1]$. Now consider a function $\Theta :(0,\infty
)\rightarrow (1,\infty )$ by $\Theta (t)=e^{t}$. Define
\[
\alpha (x,y)=\left\{
\begin{array}{cc}
1 & if~~\phi (x(t),y(t))\geq 0,t\in \lbrack 0,1], \\
0 & otherwise.%
\end{array}%
\right.
\]%
Then from (\ref{6}) with $k=\frac{1}{8},$ we obtain that
\[
\alpha (x,y)\Theta (||Fx-Fy||_{\infty })\leq \lbrack \Theta
(d(x,y))]^{k}\leq \lbrack \Theta (d(x,y)+\lambda d(y,Fx))]^{k}.
\]%
Therefore the mapping $F$ is multivalued almost $\Theta$-contraction.

From (ii) there exists $x_{0}\in C[0,1]$ such that $\alpha(x_{0},Fx_{0})\geq 1$. Next, for any $x,y\in C[0,1]$ with $\alpha (x,y)\geq
1,$ we have
\begin{eqnarray*}
&&\phi (x(t),y(t))\geq 0~~~~\text{for all}~~t\in \lbrack 0,1] \\
&\Rightarrow &\phi (Fx(t), Fy(t))\geq 0~~~~\text{for all}~~t\in \lbrack 0,1] \\
&\Rightarrow &\alpha (Fx, Fy)\geq 1,
\end{eqnarray*}%
and hence $F$ is $\alpha $-admissible. It follows from Theorem \ref{Th5.3}
that $F$ has a fixed point $x$ in $C([0,1])$ which in turns is the solution
of (\ref{2}).
\end{proof}
\end{theorem}
\section{Conclusion}
This paper is concerned with the existence and uniqueness of the best proximity point results for Brinde type contractive conditions via auxiliary function $\Theta\in\Omega$ in the framework of complete metric spaces. Also, some fixed point results as a special cases of our best proximity point results in the relevant contractive conditions are studied. Moreover, the corresponding fixed point results are obtained. An example is discussed to show the significance of the investigation of this paper. An application to a nonlinear differential equation is presented to illustrate the usability of the new theory.

\subsection*{Acknowledgments}
This paper was funded by the University of Sargodha, Sargodha funded research project No. UOS/ORIC/2016/54. The first and third author,
therefore, acknowledges with thanks UOS for financial support.

\end{document}